\documentclass[12pt]{article}

\usepackage{amsmath, epsfig, cite}
\usepackage{amssymb}
\usepackage{amsfonts}
\usepackage{latexsym}

\newtheorem{thm}{Theorem}[section]

\newcommand{\pf}{\noindent{\it Proof.} }

\numberwithin{equation}{section}

\begin{document}

\nocite{*}
\begin{center}
{\Large\bf Proof of a recent conjecture of Z.-W. Sun}
\end{center}

\vskip 2mm \centerline{Song Guo and Victor J. W. Guo\footnote{Corresponding author.}}
\begin{center}
{School of Mathematical Sciences, Huaiyin Normal University, Huai'an, Jiangsu 223300,
 People's Republic of China\\
{\tt guosong77@hytc.edu.cn,  jwguo@hytc.edu.cn} }

\end{center}


\vskip 0.7cm \noindent{\bf Abstract.} The polynomials $d_n(x)$ are defined by
\begin{align*}
d_n(x)&=\sum_{k=0}^n{n\choose k}{x\choose k}2^k.
\end{align*}
We prove that, for any prime $p$, the following congruences hold modulo $p$:
\begin{align*}
\sum_{k=0}^{p-1}\frac{{2k\choose k}}{4^k} d_k\left(-\frac{1}{4}\right)^2
&\equiv
\begin{cases}
2(-1)^{\frac{p-1}{4}}x,&\text{if $p=x^2+y^2$ with $x\equiv 1\pmod{4}$,}\\
0,&\text{if $p\equiv 3\pmod{4}$,}
\end{cases} \\[5pt]
\sum_{k=0}^{p-1}\frac{{2k\choose k}}{4^k} d_k\left(-\frac{1}{6}\right)^2
&\equiv 0, \quad\text{if $p>3$,} \\[5pt]
\sum_{k=0}^{p-1}\frac{{2k\choose k}}{4^k} d_k\left(\frac{1}{4}\right)^2
&\equiv
\begin{cases}
0,&\text{if $p\equiv 1\pmod{4}$,}\\
(-1)^{\frac{p+1}{4}}{\frac{p-1}{2}\choose \frac{p-3}{4}},&\text{if $p\equiv 3\pmod{4}$.}
\end{cases} \\
\sum_{k=0}^{p-1}\frac{{2k\choose k}}{4^k} d_k\left(\frac{1}{6}\right)^2
&\equiv 0, \quad\text{if $p>5$.}
\end{align*}
The $p\equiv 3\pmod{4}$ case of the first one confirms a conjecture of Z.-W. Sun, while the second
one confirms a special case of another conjecture of Z.-W. Sun.

\vskip 3mm \noindent {\it Keywords}: Delannoy number; congurence; Fermat's little theorem

\vskip 2mm
\noindent{\it MR Subject Classifications}: 11A07, 11B65, 05A10

\section{Introduction}
It is well known that the {\it Delannoy number}
\begin{align*}
\sum_{k=0}^{n}{n\choose k}{m\choose k}2^k=\sum_{k=0}^{n}{n\choose k}{n+m-k\choose n} 
\end{align*}
counts lattice paths from $(0,0)$ to $(m,n)$ using only single steps east $(1,0)$, north $(0,1)$, or northeast $(1,1)$.
Recently, Z.-W. Sun \cite{Sun} introduced the following polynomials
\begin{align*}
d_n(x)&=\sum_{k=0}^n{n\choose k}{x\choose k}2^k.
\end{align*}
and established some interesting supercongruences involving $d_n(x)$, such as
\begin{align*}
\sum_{k=0}^{p-1}(-1)^k d_k(x)^2 \equiv (-1)^{\langle x\rangle_p} \pmod{p^2},
\end{align*}
where $p$ is an odd prime and $\langle x\rangle_p$ denotes the least non-negative integer $r$ with $r\equiv x\pmod{p}$.
He  also made several interesting conjectures on congruences involving $d_n(x)$, such as (see \cite[Conecture 6.2]{Sun})
\begin{align}
\sum_{k=0}^{p-1}\frac{{2k\choose k}}{4^k} d_k\left(-\frac{1}{6}\right)^2
\equiv \frac{p}{3}\left(\frac{p}{3}\right)\left(4\left(\frac{-2}{p}\right)-1\right)\pmod{p^2}, \label{eq:zero}
\end{align}
where $p>3$ is a prime and $(\frac{\cdot}{p})$ is the Legendre symbol.

In this paper, we shall prove the following result.
\begin{thm}\label{thm:main}
Let $p$ be an odd prime. Then modulo $p$,
\begin{align}
\sum_{k=0}^{p-1}\frac{{2k\choose k}}{4^k} d_k\left(-\frac{1}{4}\right)^2
&\equiv
\begin{cases}
2(-1)^{\frac{p-1}{4}}x,&\text{if $p=x^2+y^2$ with $x\equiv 1\pmod{4}$,}\\
0,&\text{if $p\equiv 3\pmod{4}$,}
\end{cases} \label{eq:main-1} \\[5pt]
\sum_{k=0}^{p-1}\frac{{2k\choose k}}{4^k} d_k\left(-\frac{1}{6}\right)^2
&\equiv 0, \quad\text{if $p>3$,} \label{eq:main-2} \\[5pt]
\sum_{k=0}^{p-1}\frac{{2k\choose k}}{4^k} d_k\left(\frac{1}{4}\right)^2
&\equiv
\begin{cases}
0,&\text{if $p\equiv 1\pmod{4}$,}\\
(-1)^{\frac{p+1}{4}}{\frac{p-1}{2}\choose \frac{p-3}{4}},&\text{if $p\equiv 3\pmod{4}$.}
\end{cases} \label{eq:main-3} \\
\sum_{k=0}^{p-1}\frac{{2k\choose k}}{4^k} d_k\left(\frac{1}{6}\right)^2
&\equiv 0, \quad\text{if $p>5$.}  \label{eq:main-4}
\end{align}
\end{thm}

The $p\equiv 3$ case of \eqref{eq:main-1} was originally conjectured by Z.-W. Sun (see \cite[Conjecture 6.3]{Sun}), and the congruence
\eqref{eq:main-2} confirms the congruence \eqref{eq:zero} modulo $p$.

\section{Proof of Theorem \ref{thm:main}}
In a previous paper, Guo \cite[Lemma 3.1]{Guo} gives the following identity:
\begin{align*}
d_n(x)^2=\sum_{k=0}^{n}{n+k\choose 2k}{x\choose k}{x+k\choose k}4^k, 
\end{align*}
which is a special case of \cite[p.~80, (2.5.32)]{Slater} by noticing \cite[p.~31, (1.7.1.3)]{Slater}
(pointed out by Wadim Zudilin).

For any odd prime $p$ and $0\leqslant k\leqslant p-1$, it is easy to see that
$$
\frac{{2k\choose k}}{4^k}\equiv
\begin{cases}
\displaystyle(-1)^k {\frac{p-1}{2}\choose k} \pmod{p},&\text{if $0\leqslant k\leqslant\frac{p-1}{2}$,}\\[10pt]
0\pmod{p},&\text{if $\frac{p+1}{2}\leqslant k\leqslant p-1$.}
\end{cases}
$$
It follows that, for any $p$-adic integer $x$,
\begin{align}
\sum_{k=0}^{p-1}\frac{{2k\choose k}}{4^k} d_k(x)^2
&\equiv \sum_{k=0}^{\frac{p-1}{2}}(-1)^k {\frac{p-1}{2}\choose k} \sum_{j=0}^{k}{k+j\choose 2j}{x\choose j}{x+j\choose j}4^j \notag\\
&=(-1)^{\frac{p-1}{2}}\sum_{j=0}^{\frac{p-1}{2}}{x\choose j}{x+j\choose j}{j\choose \frac{p-1}{2}-j}4^j \pmod{p} \label{eq:dsqure-x}
\end{align}
by noticing the Chu-Vandermonde identity
\begin{align*}
\sum_{k=j}^{\frac{p-1}{2}}(-1)^k {\frac{p-1}{2}\choose k}{k+j\choose 2j}
=(-1)^{\frac{p-1}{2}}{j\choose \frac{p-1}{2}-j}.
\end{align*}

It is not difficult to see that
\begin{align*}
{-\frac{1}{4}\choose j}{-\frac{1}{4}+j\choose j}&=(-1)^j\frac{{4j\choose 2j}{2j\choose j}}{64^j}
\equiv 0\quad\text{for}\ \frac{p}{4}\leqslant j\leqslant\frac{p-1}{2},\\[5pt]
{-\frac{1}{6}\choose j}{-\frac{1}{6}+j\choose j}&=(-1)^j\frac{{6j\choose 3j}{3j\choose j}}{432^j}
\equiv 0\quad\text{for}\ \frac{p}{6}\leqslant j\leqslant\frac{p-1}{2},\\[5pt]
{\frac{1}{4}\choose j}{\frac{1}{4}+j\choose j}&=(-1)^{j-1}\frac{(4j+1){4j\choose 2j}{2j\choose j}}{(4j-1)64^j}
\equiv 0\quad\text{for}\ j=\frac{p-1}{4}\ \text{or}\ \frac{p+3}{4}\leqslant j\leqslant\frac{p-1}{2},\\[5pt]
{\frac{1}{6}\choose j}{\frac{1}{6}+j\choose j}&=(-1)^{j-1}\frac{(6j+1){6j\choose 3j}{3j\choose j}}{(6j-1)432^j}
\equiv 0\quad\text{for}\ \frac{p+3}{6}\leqslant j\leqslant\frac{p-1}{2},\end{align*}
and ${j\choose \frac{p-1}{2}-j}=0$ for $0\leqslant j<\frac{p-1}{4}$. Letting $x=\mp\frac{1}{6}$ in \eqref{eq:dsqure-x}, we immediately obtain
\eqref{eq:main-2} and \eqref{eq:main-4}.
Letting $x=\mp\frac{1}{4}$ in \eqref{eq:dsqure-x}, we get
\begin{align*}
\sum_{k=0}^{p-1}\frac{{2k\choose k}}{4^k} d_k\left(-\frac{1}{4}\right)^2
&\equiv
\begin{cases}
(-1)^{\frac{p-1}{4}}\frac{{p-1\choose \frac{p-1}{2}}{\frac{p-1}{2}\choose \frac{p-1}{4}}}{16^{\frac{p-1}{4}}}\pmod{p},&\text{if $p\equiv 1\pmod{4}$,}\\
0\pmod{p},&\text{if $p\equiv 3\pmod{4}$,}\\[5pt]
\end{cases} \\[10pt]
\sum_{k=0}^{p-1}\frac{{2k\choose k}}{4^k} d_k\left(\frac{1}{4}\right)^2
&\equiv
\begin{cases}
0\pmod{p},&\text{if $p\equiv 1\pmod{4}$,}\\
(-1)^{\frac{p+1}{4}}\frac{(4p+2){p+1\choose \frac{p+1}{2}}{\frac{p+1}{2}\choose \frac{p+1}{4}}}{4p\cdot16^{\frac{p+1}{4}}}\pmod{p},&\text{if $p\equiv 3\pmod{4}$,}\\[5pt]
\end{cases}
\end{align*}
We now suppose that $p$ is a prime such that $p\equiv 1\pmod{4}$  and $p=x^2+y^2$ with $x\equiv 1\pmod{4}$. Then
by the Beukers-Chowla-Dwork-Evans congruence \cite{CDE,Pan2}:
\begin{align*}
{\frac{p-1}{2}\choose \frac{p-1}{4}}\equiv\frac{2^{p-1}+1}{2}\left(2x-\frac{p}{2x}\right) \pmod{p^2},
\end{align*}
and Fermat's little theorem, we have ${\frac{p-1}{2}\choose \frac{p-1}{4}}\equiv 2x\pmod{p}$. Moreover, we have
$$
{p-1\choose \frac{p-1}{2}}\equiv 16^{\frac{p-1}{4}}\equiv 1\pmod{p}.
$$
This proves \eqref{eq:main-1}. Finally, suppose that $p$ is a prime with $p\equiv 3\pmod{4}$. Then
$$
\frac{(4p+2){p+1\choose \frac{p+1}{2}}{\frac{p+1}{2}\choose \frac{p+1}{4}}}{4p\cdot16^{\frac{p+1}{4}}}
\equiv {\frac{p-1}{2}\choose \frac{p-3}{4}}\pmod{p}.
$$ This proves \eqref{eq:main-3}.

\vskip 5mm \noindent{\bf Acknowledgment.} This work was partially supported by the National Natural Science Foundation of China (grant no. 11371144).

\end{document}